\documentclass[12pt]{article}
\usepackage{amsmath,amssymb,amstext,dsfont,fancyvrb,float,fontenc,graphicx,subfigure,theorem,hyperref}
\usepackage[hyphenbreaks]{breakurl}

\usepackage[utf8]{inputenc}

\usepackage[letterpaper]{geometry}
\title{\bf The Integral of Secant and Stereographic Projections of Conic Sections}
\author{George Jennings, David Ni, Wai Yan Pong, Serban Raianu\\\\
California State University, Dominguez Hills\\1000 E Victoria St, Carson, CA 90747}

\begin{document}
\setcounter{page}{1}
\maketitle

\begin{abstract}
\noindent We show that the four best-known substitutions used to integrate secant rely on different rational parametrizations of conic sections, coming from stereographic projections. We also investigate other possible explanations for the four substitutions, as well as connections between the buildups leading to them.
\end{abstract}
\section{Introduction}
 The integral of secant, which competes for the unofficial title of \textit{The Most Interesting Integral in the World of Elementary Functions}, can be taught in several ways \cite{wiki,mit,clp,hardy,ste}. 

We start by briefly reviewing some history \cite{rt}. 
The secant was integrated numerically in 1569 by Gerardus Mercator \cite{merc} to construct the famous map projection that bears his name.  Mercator's projection is especially useful for navigators because it depicts angles accurately, with zero distortion.   Let $\theta = $ longitude and $\phi = $ latitude of a point on Earth. On Mercator's map the $x$-coordinate of the corresponding point on the map is the longitude, and its $y$-coordinate depends only on the latitude  
  $$x=\theta, \quad y=y(\phi).$$
Consider an infinitesimal square on the earth with vertices at $(\theta,\phi)$, $(\theta+d\theta,\phi)$, $(\theta,\phi+d\phi)$, and $(\theta+d\theta,\phi+d\phi)$.  Its sides have equal length, 
$$\text{side}=ds=r\cos\phi\;d\theta = r\;d\phi$$
where $r$ is the Earth's radius.  To avoid distorting angles the corresponding points on the map must also form an infinitesimal square.  Its horizontal measure is 
$$dx=d\theta = \frac{ds}{r\cos\phi}$$
so its vertical measure must be the same,
$$dy=\frac{ds}{r\cos\phi}=\frac{rd\phi}{r\cos\phi}=\sec\phi\;d\phi.$$
Adding up these infinitesimal distances $dy$ on the map as one moves away from the equator, the $y$-coordinate of a point at latitude $L$ will be a function of $L$, say $y=f(L)$. Then we have 
$$f(L) = \int_0^{f(L)} dy = \int_0^{L} \sec\phi\;d\phi.$$
Taking derivatives shows that $f'(L)=\sec L$. Note that this gives us a practical reason to integrate the secant function but it does not help us find $f(L)$. This is done by one of the following four available substitutions described in the next section.
\section{The Four Substitutions Used to Integrate Secant}
\begin{enumerate}
\item{The Gregory substitution  \cite{greg,inglis}, $u=\sec x+\tan x$. This is the method used in most modern calculus textbooks \cite{ste,clp}, and is usually justified by a mysterious trick: multiplying and dividing $\sec x$ by $\sec x +\tan x$. Here is an explanation (found independently many times, see \cite{mit,sch,strauss}) for how to come up with this: integration is about spotting derivatives. Not being able to think of a function whose derivative is secant, we go for the next best thing, i.e. finding a function which at least has a secant in its derivative. It turns out there are only two of them:
\begin{eqnarray*}
(\tan x)' &=&\sec^2 x\\
(\sec x)'&=&\sec x\tan x
\end{eqnarray*}
Seeing them nicely lined up like this prompts us to add the two equalities, and by doing so we get:
$(\sec x+\tan x)'=\sec x(\sec x+\tan x)$, so 
\begin{equation}\label{intsec}
\int\sec x\,\mathrm{d}x=\ln |\sec x+\tan x|+C.
\end{equation}
This was a conjecture from 1645, based on information from logarithmic tables (and known to many people, including Newton), and it was solved by James Gregory 23 years later \cite{rt}. As explained by Turnbull in \cite{tur}, Gregory actually proves (using only geometry), that
$$\int_0^{\alpha}\sec\theta\;\;{\rm d}\theta=\int_{1/\sqrt{2}}^{(\sec\alpha+\tan\alpha)/\sqrt{2}}\;\;\frac{{\rm d}u}{u}=\log(\sec\alpha+\tan\alpha),$$
and this is how he arrives at the substitution.}
\item{The substitution used by Newton's teacher, Barrow \cite{barr}, $u=\sin x$. This is set up by multiplying and dividing $\frac{1}{\cos x}$ by $\cos x$, which, as we can see from the explanation of the Gregory substitution in \cite{tur}, was the first step in Gregory's chain of simple and double integrals. }
\item{The Weierstrass \cite{weier} substitution, $t=\tan\left(\frac{x}{2}\right)$ \cite{joh}.}
\item{The modified Weierstrass substitution \cite{hardy}, $s=\tan\left(\frac{x}{2}+\frac{\pi}{4}\right)$, which leads to the integral of $\frac{1}{s}$.}
\end{enumerate}
The Barrow and Weierstrass substitutions both end up by using the method of partial fractions to integrate $\frac{1}{1-u^2}$ (or $\frac{1}{1-t^2}$, respectively).
It is believed that Barrow was the first to use partial fractions in integration \cite{rt}. In the North American curriculum partial fractions are usually taught after trigonometric integrals \cite{ste}, but the following well known trick helps us avoid partial fractions in this particular case, so we can teach the integral of secant using any of the four methods at any time. Moreover, the antiderivative is the same as the one obtained when using partial fractions.
\begin{align*}
\int \frac{{\rm d}u}{1 - u^2} & = \int\frac{(1-u+u){\rm d}u}{(1+u)(1-u)} = \dfrac12\int \left(\frac{2}{1+u} - \frac{-2u}{1-u^2}\right)\,{\rm d}u \\	
& = \frac12 \left( 2 \ln\left|1 + u\right|-\ln \left|1-u^2\right| \right) + C = \frac12 \ln\left|\frac{1+u}{1-u}\right| + C	
\end{align*}
A similar but sneakier trick produces the partial fraction decomposition by splitting the 1 in the numerator into $1/2+1/2$, then adding and subtracting $u/2$ and separating the fractions. 

We would to recall here that the method of partial fractions is a spectacular and surprising application of  the arithmetic of the ring of polynomials in one indeterminate with real or complex coefficients to calculus \cite[p. 165]{clp}. 

\section{Connections Between the Four Substitutions and Other Ideas}
We now discuss another idea for integrating secant and a different buildup to the Gregory substitution, both coming from Liouville's Theorem on integration in finite terms, as proved in \cite{mrosen}. (By the way, the proof makes frequent use of partial fractions.) This is a theorem in differential algebra which explains why functions like $e^{x^2}$ or $\sin x^2$ do not have antiderivatives which can be expressed as compositions of elementary functions (which are compositions of the functions we study in calculus). We give the full statement of the theorem, but keep in mind that an easy example of a differential field is the field of rational fractions, and for a rational fraction $u$ we can think of $u'$ as the derivative of $u$ regarded as a rational function.\\[1mm]

\noindent {\bf Liouville's Theorem.} Let $F$ be a differential field of characteristic zero and $\alpha\in F$. If the
equation $y'=\alpha$ has a solution in some elementary differential extension field of $F$
having the same subfield of constants, then there are constants $c_1,\ldots ,c_n\in F$ and
elements $u_1,\ldots,u_n,v\in F$ such that
$$\alpha=\sum_{i=1}^n c_i\frac{u_i'}{u_i}+v'.$$

Now, if $\sec x$  has an elementary indefinite integral, it would be of the form in Liouville's Theorem, and so, since we can't think of a function whose derivative is secant, guessing $u(x)$ in the simplest case $$\sec x = \frac{u'(x)}{u(x)} =\frac{{\rm d}}{{\rm d}x} \ln |u(x)|$$  would be a reasonable thing to try first. We therefore try to find $u(x)$ if $u'(x)\cos x=u(x)$. Since we are looking for $u'(x)\cos x$, it makes sense to differentiate $h(x)=u(x)\cos x$.
\begin{eqnarray*}
h'(x) &=& -u(x)\sin x +u'(x)\cos x\\
&=& -u(x)\sin x +u(x)\\
&=& -\frac{h(x)}{\cos x}\sin x +\frac{h(x)}{\cos x}\\
&=& h(x)\frac{1-\sin x}{\cos x}\\
&=& h(x)\frac{\cos x}{1+\sin x}.
\end{eqnarray*}
This gives
$$\ln |h(x)|=\ln (1+\sin x),$$
and so 
\begin{equation*}
\ln |u(x)|=\ln \left|\frac{h(x)}{\cos x}\right|=\ln |\sec x+\tan x|,
\end{equation*}
which means taking $u(x)=\sec x+\tan x$ works.

This also tells us that the buildup to Barrow's substitution, combined with the avoidance of partial fractions described above, can also lead  to Gregory's substitution. Keep in mind that Barrow's buildup involved multiplying the top and the bottom by $\cos x$, then adding and subtracting  $\sin x$ in the numerator. If instead we reverse the order of these two actions, i.e. we start by adding and subtracting $\sin x$ on top, then we get:
\begin{equation}\label{one}
\frac{1}{\cos x}=\frac{\sin x}{\cos x}+\frac{1-\sin x}{\cos x}=\frac{\sin x}{\cos x}+\frac{\cos x}{1+\sin x}
\end{equation}
(multiplying and dividing by $\cos x$ is now hidden in the identity $$\frac{1-\sin x}{\cos x}=\frac{\cos x}{1+\sin x},$$ this may be seen after dividing both of its sides by $1-\sin x$),
and (\ref{one})  immediately gives (\ref{intsec}).

Now we can get yet another buildup leading to Gregory's substitution. Start with (\ref{one}), and keep going instead of integrating. We get
\begin{eqnarray*}
\frac{1}{\cos x} &=& \frac{\sin x}{\cos x}+\frac{\cos x}{1+\sin x}\\
&=& \tan x+\frac{1}{\sec x+\tan x}\\
&=& \frac{\tan x(\sec x+\tan x)+(\sec^2 x-\tan^2 x)}{\sec x+\tan x}\\
&=&  \frac{(\sec x +\tan x)\sec x}{\sec x+\tan x}.
\end{eqnarray*}

We end this section by noting that we can also derive (\ref{one}) by adding and subtracting $\tan x$ to $\sec x$ instead of $\sin x$ to 1 on top of $\frac{1}{\cos x}$  (this can be also independently motivated by the identity $\sec^2 x-\tan^2 x=1$). 
\begin{equation}\label{two}
\sec x=\tan x+\sec x-\tan x=\tan x+\frac{1}{\sec x+\tan x}=\frac{\sin x}{\cos x}+\frac{\cos x}{1+\sin x}
\end{equation}
\section{The Substitutions Used in the Integral of Secant and Stereographic Projections}
The {\em stereographic projection} of a curve from a point  $O$ to a line is the function that maps any point $P$ on the curve  to the intersection of that line with the secant $OP$ to the curve. The usual projection parallel to a line is a particular case (the point in this case is a point at infinity).

The miraculous fact that the Weierstrass substitution always works, and turns any rational fraction in trigonometric functions into an integrable rational function in $t$ is a direct consequence of the fact that a circle is a curve of genus 0 (it has no ``handles"), which also also means it has a rational parametrization. The Weierstrass substitution 
$t=\tan(\theta/2)$ arises naturally from the stereographic projection of the unit circle from $(-1,0)$ to the $y$-axis.  Let 
$$(x,y)=(\cos\theta, \sin\theta)\neq (-1,0)$$
be a point on the unit circle with center at the origin.  The line connecting $(x,y)$ to $(-1,0)$ intersects the $y$-axis at a point $(0,t)$.  Elementary geometry (the triangle with vertices $(x,y),(-1,0),(0,0)$ is isosceles)  says the angle between this line and the $x$-axis is $\theta/2$.  {\em Stereographic projection} of the circle from the point $(-1,0)$ to the $y$-axis is the function that maps $(x,y)$ to $(0,t)$.  

\begin{figure*}[h]
  \begin{center}
  \includegraphics[scale=.5]{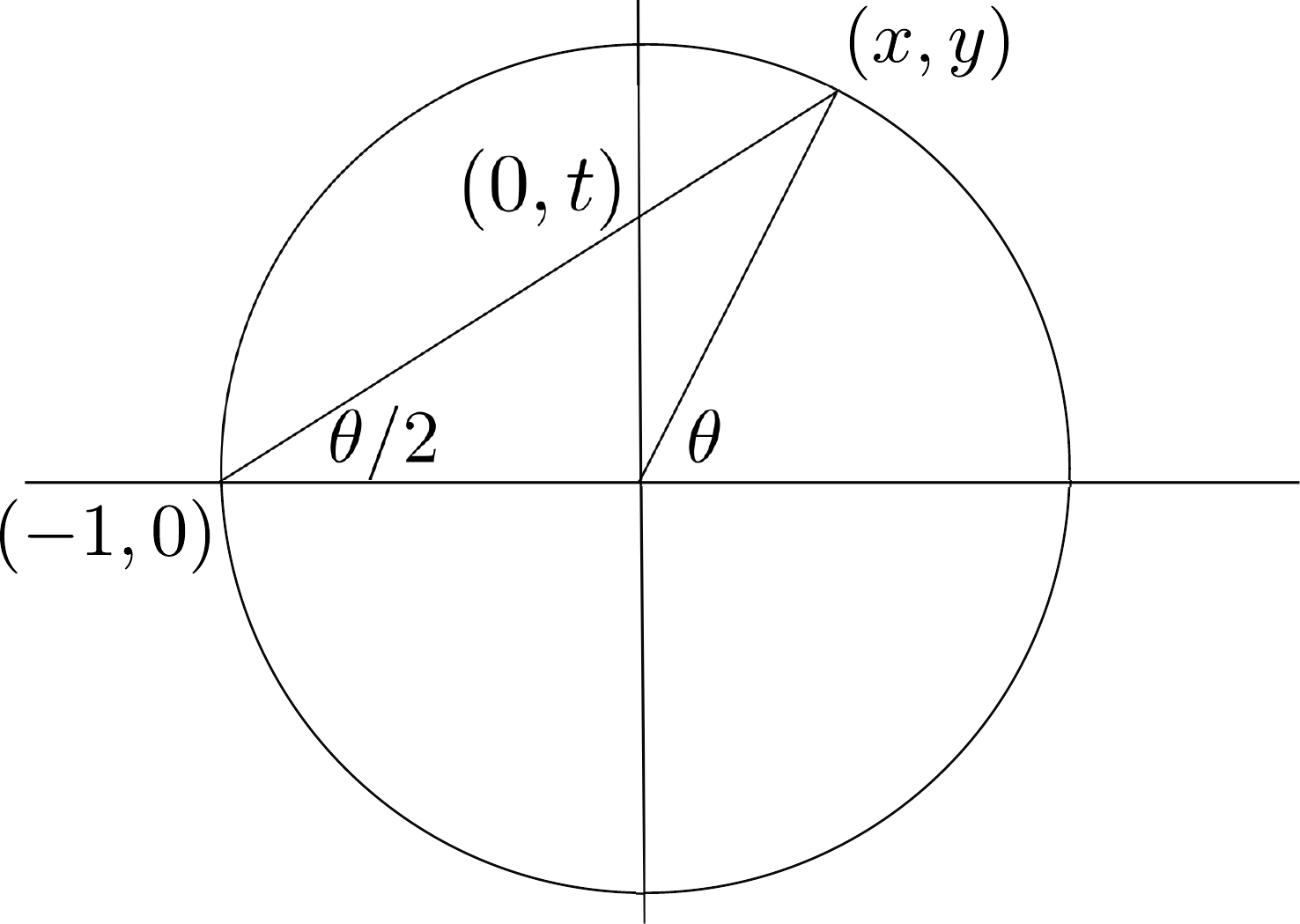}
\end{center}
\caption{\label{fig1} A rational parametrization of the circle}
\end{figure*}

The slope of the line is $t=y/(x+1)=\tan(\theta/2)$.  Plug $y=t(x+1)$ into the equation $x^2+y^2-1=0$, divide by $x+1$, and solve for $x$ to obtain
\begin{equation}\label{stereo}
(x,y)=(\cos\theta,\sin\theta)=\left( \frac{1-t^2}{1+t^2},\frac{2t}{1+t^2} \right),
\end{equation} 
a \emph{rational} parametrization of the unit circle.  Differentiating
$$y=\sin\theta = \frac{2t}{1+t^2}$$
and using $x=(1-t^2)/(1+t^2)$ one obtains
\begin{gather*}
 dy = \cos\theta\;d\theta = \frac{2(1-t^2)}{(1+t^2)^2}\;dt = \frac{2\cos\theta}{1+t^2}\; dt \\
 \intertext{hence}
 d\theta = \frac{2\;dt}{1+t^2}.
\end{gather*}
These substitutions help convert  \emph{any} rational trigonometric differential $ R(\cos\theta,\sin\theta) \;d\theta$, where $R(x,y)$ is a quotient of polynomials (including $\sec\theta d\theta$), to a rational differential $S(t)\;dt$ that is sometimes easier to integrate.

We also note (see \cite[p. 4]{st}) that this rational parametrization of the circle also produces a recipe for constructing Pythagorean triples of integers $X,Y,Z$ such that $X^2+Y^2=Z^2$ (also see \cite{pyth} for a biography of Pythagoras and  \cite[Theorem 6.3, p.162]{str} for a purely algebraic description of all primitive --- i.e. with no common factors --- Pythagorean triples). Let $t=a/b$ where $a$ and $b\neq 0$ are integers. Plug $t$ into the stereographic projection formula (\ref{stereo}) to obtain $(x(t),y(t))$ satisfying $x^2+y^2=1$ then clear denominators to obtain $X^2+Y^2=Z^2$ where $X,Y,Z$ are integers.  Conversely, given a Pythagorean triple $(X,Y,Z)$ one may reverse this process to find a rational $t$ that will produce $(X,Y,Z)$ or an integer multiple of it.

\medskip
The modified Weierstrass substitution $s=\tan(\theta/2+\pi/4)$ also uses a stereographic projection of the unit circle, but from the point $(0,1)$ to the $x$-axis. 
\begin{figure*}[h]
  \begin{center}
    \includegraphics[scale=.5]{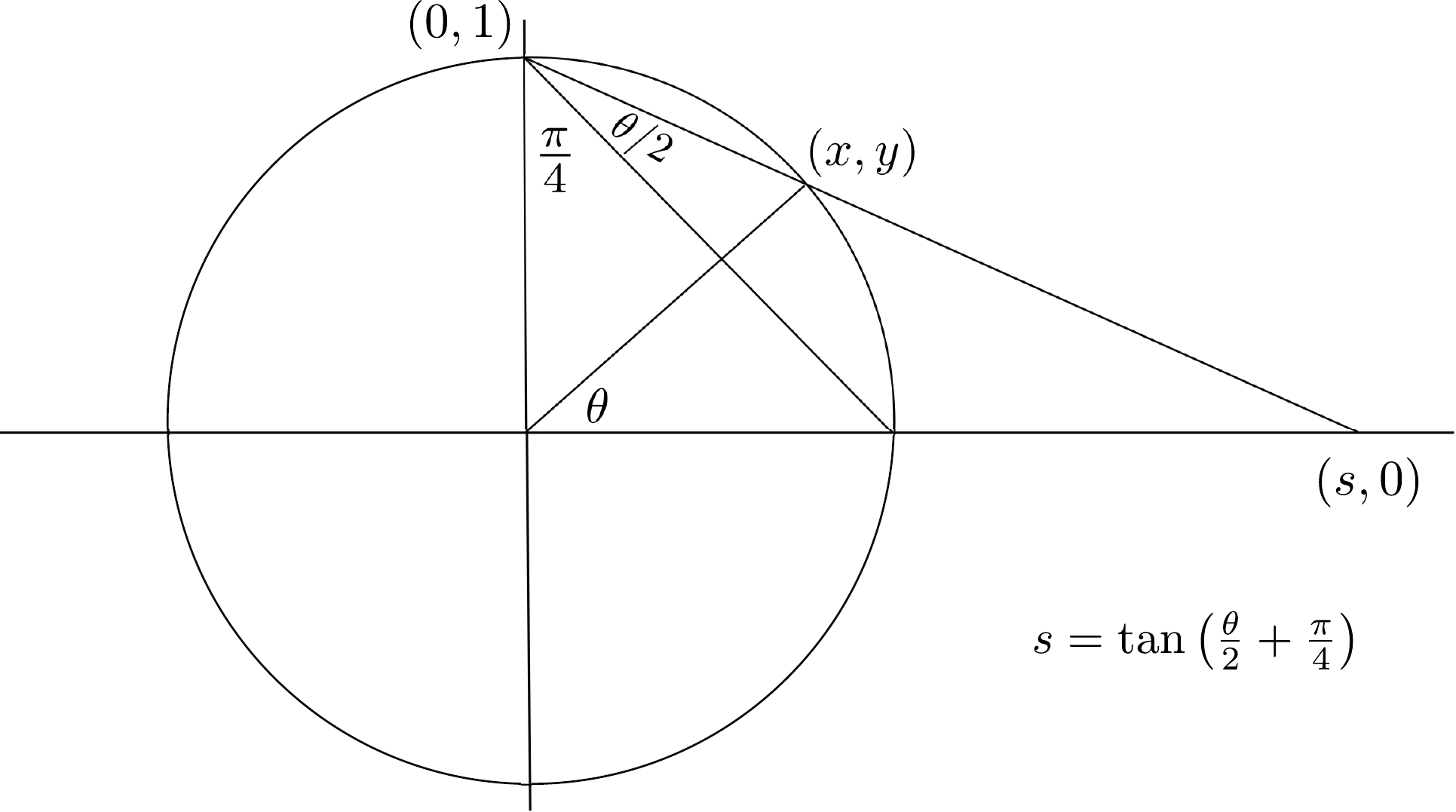}
  \end{center}
\caption{\label{fig2} Another rational parametrization of the circle}
\end{figure*}

Let $(\cos\theta,\sin\theta)=(x,y)\neq(0,1)$.  The line through $(x,y)$ and $(0,1)$ intersects the $x$-axis at $(0,s)$, where $s=\tan(\theta/2+\pi/4)$. Its slope is $-1/s= (y-1)/x$ so $x=s(1-y)$.  Plug that into the equation $x^2+y^2-1=0$, divide by $1-y$, and solve for $y$ to obtain
$$(x,y)=(\cos\theta,\sin\theta)=\left( \frac{2s}{s^2+1},\frac{s^2-1}{s^2+1} \right),$$ 
another rational parametrization of the unit circle.  As before, computing $dy=\cos\theta\; d\theta$ yields a formula for $d\theta$
$$d\theta = \frac{2}{s^2+1}\; ds.$$
These formulas make integrating 
$$\int\sec\theta\; d\theta = \int \frac{s^2+1}{2s} \frac{2}{s^2+1}\; ds = \int \frac{ds}{s} = \ln\lvert s\rvert +C $$
even simpler, without partial fractions.  Finally, using $x^2+y^2=1$,  
\begin{equation*}
s = \frac{x}{1-y} = \frac{x(1+y)}{1-y^2}=\frac{x(1+y)}{x^2}=\frac{1}{x}+\frac{y}{x}=\sec\theta + \tan\theta
\end{equation*}
produces the formula that appears in the back of your textbook. This shows that the Gregory substitution is basically the same as the modified Weierstrass. Nevertheless, we give a separate geometric argument for the Gregory substitution:

The parametrization of the right branch of the hyperbola $x^2-y^2=1$ is $(x,y)=(\sec\theta,\tan\theta)$.
The line at infinity $z=0$ intersects the homogenized hyperbola $x^2-y^2=z^2$ at the two points at infinity $P_-=[1:1:0]$ and $P_+=[1:-1:0]$. 
The stereographic projection of the hyperbola from $P_+$ to the $x$-axis (which is the projection down a line  $x+y=u$ parallel to the asymptote $x+y=0$) sends the point $(x,y)$ on the hyperbola to the point $(x+y,0)=(\sec\theta+\tan\theta,0)$, and we choose the parameter to be $u=x+y=\sec\theta+\tan\theta$. 
Solving for the point where the line $x+y=u$ intersects the hyperbola, we obtain  
$$(x(u),y(u))=\left(\frac{u^2+1}{2u},\frac{u^2-1}{2u}\right)=(\sec\theta,\tan\theta).$$
This rational parametrization of the hyperbola  does the job.

\begin{figure*}[h]
  \begin{center}
\includegraphics[scale=.5]{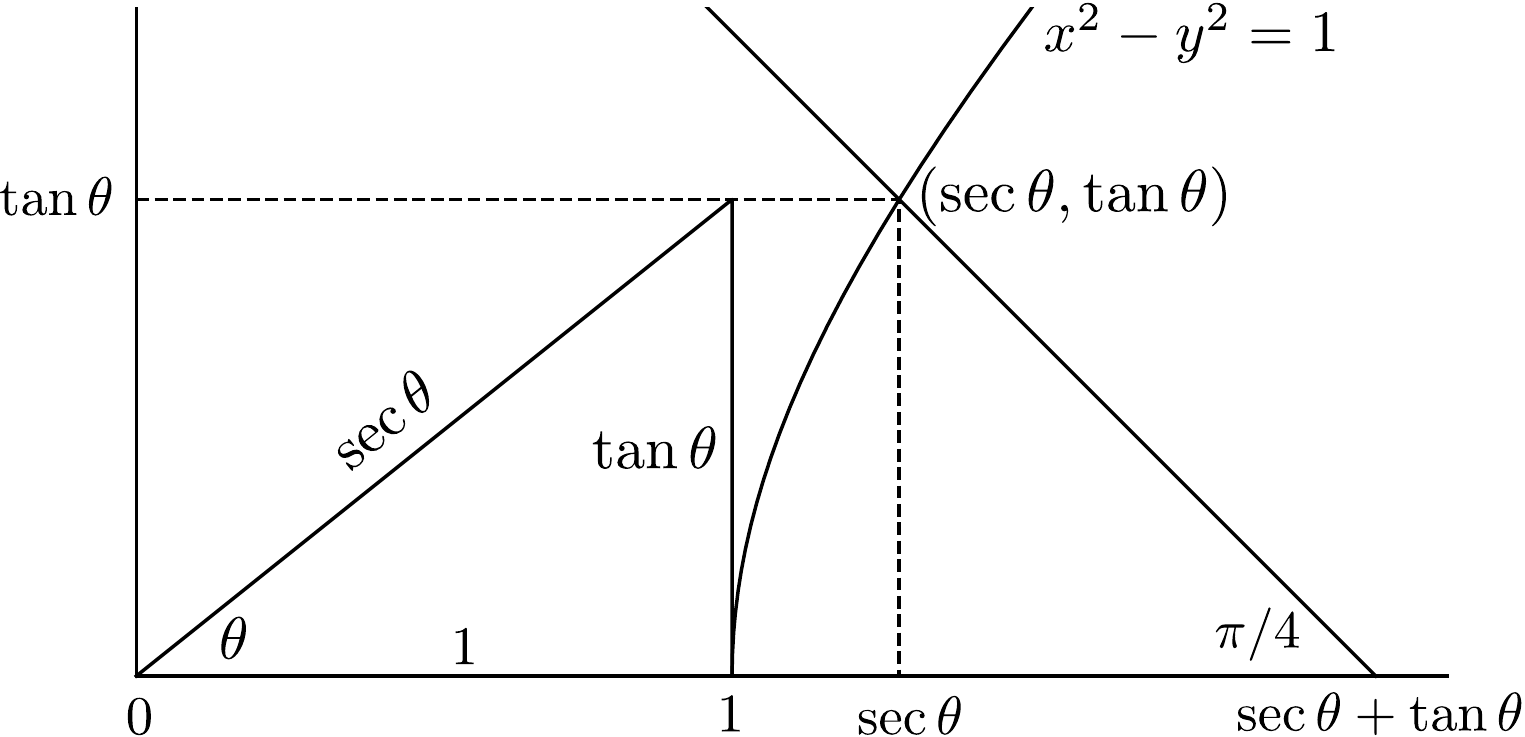}
\end{center}
\caption{\label{fig3} A rational parametrization of the hyperbola}
\end{figure*}

\noindent We have
$$\mathrm{d} u=(\sec\theta\tan\theta+\sec^2\theta)\mathrm{d}\theta=u\sec\theta\mathrm{d}\theta,$$
so
$$\frac{\mathrm{d} u}{u}=\sec\theta\;\mathrm{d}\theta,$$
and this provides another natural buildup to Gregory's substitution.

Other rational parametrizations of the hyperbola also work.  One possible parameter is $v$ in the point $(v,v)$ which is the stereographic projection from $P_+$ of the point $(x,y)=(\sec\theta,\tan\theta)$ on the hyperbola to the asymptote $y=x$.  Another is $w=\sqrt{2v}$, which is the distance from $(v,v)$ to the origin.  Both of these two parameters work just as well as the $u$ above, because $v=u/2$, and $w=u/\sqrt{2}$. Other choices include the parameters we have used before: $t=\tan(\theta/2)$ and $s=\tan(\theta/2+\pi/4)$, but these end up requiring partial fractions.  The reader may wish to work out formulas for parametrizing the hyperbola with these parameters as an exercise. 

We now see that all methods for integrating the secant function have something in common: they rely on certain rational parametrizations of the circle or the hyperbola. The figure below puts all these parametrizations together on the same picture. This shows how our previous three pictures are related.

\begin{figure*}[h]
  \begin{center}
\includegraphics[scale=0.8]{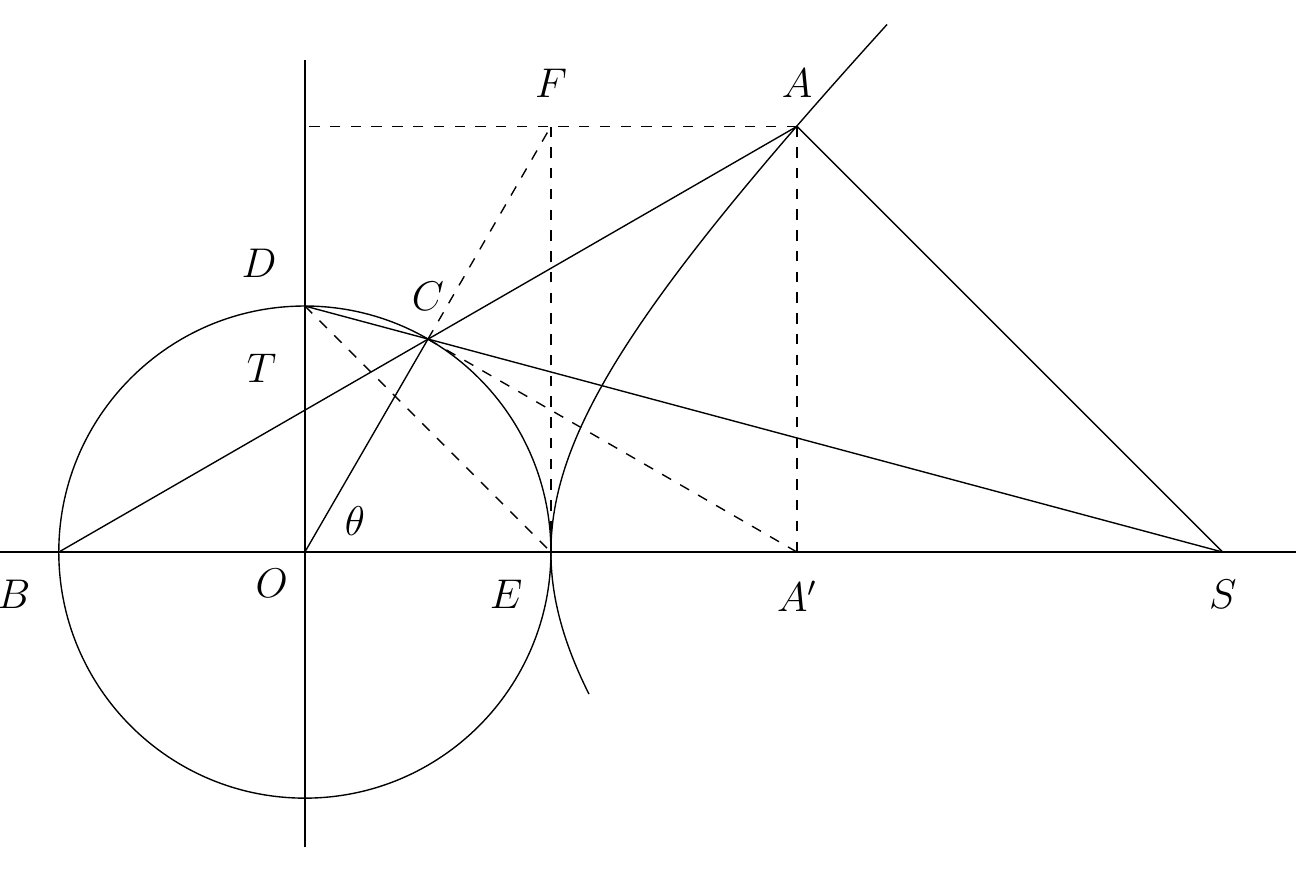}
\end{center}
\caption{\label{fig4} All previous three rational parametrizations}
\end{figure*}

The relevant points are $B(-1,0)$, $D(0,1)$, $C(\cos\theta,\sin\theta)$, $A(\sec\theta,\tan\theta)$, and the
parameters are:
\begin{itemize}
\item $OT=t=\tan(\theta/2)$ (Weierstrass)
\item $OS=s=\tan(\theta/2+\pi/4)=\sec\theta+\tan\theta$ (modified Weierstrass, Barrow --- see (\ref{one}) or (\ref{two}), Gregory)
\end{itemize}

This picture is bursting with interesting geometry problems; see for example \cite{jnpr}. The most striking of these geometric facts is that the stereographic projection of the circle from the point $D(0,1)$ to the $x$-axis and the stereographic projection of the hyperbola from the point at infinity $P_+[1:-1:0]$ to the $x$-axis (recall that the latter is the projection of the hyperbola to the $x$-axis  parallel to the asymptote $x+y=0$) turn out to be the same  if we identify the circle and the hyperbola via corresponding points ($A$ and $C$).


\section*{Acknowledgments}
We thank Frank Miles for helpful comments.

\vfill\eject

\end{document}